\documentstyle[12pt, twoside]{article}
\newcommand{\B} {\lambda}

\begin{document}
\date{}
\newtheorem{thm}{Theorem}[section]
\newtheorem{lem}{Lemma}[section]
\newtheorem{rem}{Remark}[section]
\newtheorem{prop}{Proposition}[section]
\newtheorem{cor}{Corollary}[section]
\newtheorem{defn}{Definition}[section]

\markboth{\hfill Liouville energy \hfill }{\hfill X. Chen \& M.
Zhu \hfill}

\title
{Liouville energy on a topological two sphere}\author  { XiuXiong Chen\\Department of Mathematics\\
University of Wisconsin-Madison \\Madison WI 53706-1388\\
Meijun Zhu\\
Department of Mathematics\\
The University of Oklahoma\\
Norman, OK 73019\\
}

\newcommand{\A} {\alpha}

\maketitle
\input { amssym.def}
\input { amssym.def}
\pagenumbering{arabic}
\newcommand{\bg} {\begin{equation}}
\newcommand{\ede} {\end{equation}}
\newcommand{\M} {\Omega}
\newcommand{\E}{\epsilon}
\newtheorem{definition}{Definition}[section]

\begin{abstract}
In this paper we shall give an analytic proof of the fact that the
Liouville energy on a topological two sphere is bounded from
below. Our proof does not rely on  the uniformization theorem and
the  Onofri inequality,  thus it is essentially needed in the
alternative proof of the uniformization theorem via the Calabi
flow. Such an analytic approach also sheds light on how to obtain
the boundedness for $E_1$ energy in the study of general K\"ahler
manifolds.

\end{abstract}

\section{Introduction}
\subsection {Description of the problem}
Let $(M,g)$ be a smooth Riemann surface. For any conformal new
metric $g_1=e^{u} g$, the corresponding Liouville energy is
defined by \bg  L_g(g_1)=  \int_M \ln \frac {g_1}{g} \cdot
(R_{g_1} dV_{g_1}+ R_{g} dV_{g}) \label{1-1} \ede where $R_g$ and
$R_{g_1}$ are  twice the  Gaussian curvatures $K_{g}$ and
$K_{g_1}$ with respect to metrics $g$ and $g_1$. Since $$
R_{g_1}=e^{-u}(-\Delta_g u+R_g),
$$ the Liouville energy of metric $g_1$ can also be represented by
$$ L_g(g_1) =  \int_M  (|\nabla_g u|^2+2 R_g u) dV_{g}.
$$

It was showed by Chow \cite{Cho} and Chen \cite{Che} that Ricci
flow and Calabi flow can be viewed as the gradient flow of the
Liouville energy.  Therefore if $M$ is a topological two sphere,
it is important to know in the study of Ricci and Calabi flow that
the Liouville energy is bounded from below for any metric with
fixed volume. In fact such boundedness can be used to show the
global existence of
 Calabi flow via integral estimates (see, for example,  Chen \cite{Che}),
 thus one can give an alternative proof of
 the uniformization theorem (for the hardest case: positive
curvature case) via  Calabi flow\footnote {In a recent note, Chen,
Lu and Tian \cite{CLT} clarify that the Ricci flow did yield
another proof of the uniformization theorem for spheres even
without using current integral estimates. Nevertheless our result
not only answer a conjecture of Chen \cite{Che}, but also has
broad applications in the study of  K\"ahler manifolds (see the
historic note).}. However, present argument in the literature
seems to be a tautology: the proof of the fact that the Liouville
energy on a two-sphere is bounded from below relies on the
uniformization theorem. In fact, the proof can be summarized as
follows. Let $(M, g)$ be a topological sphere. From the
uniformization theorem we know that $(M, g)$ is conformally
equivalent to the standard sphere $(S^2, g_0)$, that is: there is
a $\varphi(x)\in C^2(S^2)$ such that $g=e^{ \varphi} g_0$. For a
given metric $g_1=e^{u} g$, since $g_1=e^{u+ \varphi} g_0$ and
$R_{g_0}=2$, we have:
$$
\begin{array}{rll}
 L_g(g_1)& = \int_M(|\nabla_g u|^2+2 R_g u) dV_{g}\\
 &=L_{g_0}(g_1)-L_{g_0}(g)\\
 &=\int_{S^2}(|\nabla_{g_0} (u+\varphi)|^2+4  (u+\varphi)) dV_{g_0}-L_{g_0}(g)\\
 &\ge 16 \pi \cdot \ln \frac 1{4 \pi} \int_{S^2} e^{(u+\varphi)}dV_{g_0}-L_{g_0}(g).
 \end{array}
$$
The last inequality follows from the well-known Onofri inequality
on $S^2$. Thus the finiteness of the volume of the metric $g_1$
implies that  the Liouville energy of such metric is bounded from
below.

\subsection{Historic note}
It is always natural to ask if one can use variational approach to
attack the existence of constant scalar curvature metrics in
certain conformal class in 2-sphere. The elliptic approach might
start with Berger \cite{B}, which eventually leads to the famous
prescribing curvature problem on $S^2$ (the Nirenberg problem).
The parabolic approach  was first studied by Hamilton \cite{H} and
Chow \cite{Cho} via Ricci flow,  and later details of the approach
via the Calabi flow were illustrated
 by the first author. The
first step of such an approach would be to find a direct analytic
proof of the fact the Liouville energy is bounded from below
first. This question was raised in \cite{Che} as a key step in
deriving a new proof for uniformization theorem in topological two
sphere via the Calabi flow.  More importantly, this question is
closely related to a far reaching program  in the study of general
K\"ahler manifolds: Can one  use variational approach to derive
the existence of constant scalar curvature  metrics in each
K\"ahler class?

     The search of constant scalar curvature metric in each K\"ahler class was suggested by
 E. Calabi  \cite{calabi82} in 1982.  This is one of the key problems in K\"ahler geometry
 and  there is extensive research in this and related problems.
 The analytic obstruction to the existence of such metric was first observed in the important work
 of Mabuchi-Bando \cite{Ba} in 1987, where they showed that the existence of K\"ahler Einstein
 metric implies the existence of lower bound of the K energy in each
K\"ahler class\footnote{On $S^2$, a K\"ahler class is equivalent a
pointwise conformal class of Riemannian metric with fixed total
area.}.  This obstruction was extended by the work of the first
author in 1998 to the case of constant scalar curvature metric  in
K\"ahler manifold with non-positive first Chern class.  More
recently, this problem is completely   settled for all K\"ahler
classes (cf. \cite{CT1} for more references).   In \cite{CT2},
Chen and Tian introduced a set of of $n+1$ energy functionals
$E_0, E_1, \cdots E_n$  in any $n-$ dimensional K\"ahler manifold
where $E_0$ is the K energy functional and $E_1$ is the Liouville
energy on $S^2.\;$  A natural question asked in \cite{Che}
(Question 7.3) and \cite{Che1}
  is whether the existence of a K\"ahler Einstein metric implies the lower bound of $E_1.$
Concentrated effort by various geometers and rapid progress are
made in this problem.  For instance,   the first author observed
that, along the K\"ahler Ricci flow on K\"ahler Einstein manifold,
this functional is bounded below.   This is answered affirmatively
by Song and Weikeov \cite{SW}  that lower bound of $E_1$ is
obstruction to the existence of K\"ahler Einstein metrics. Shortly
after, Pali \cite{Pa}   gives a new proof which considerably
simplified the Song-Weikeov's argument.  See also \cite{Va} for
other related results.  All of these recent results highlight the
importance of obtaining an apriori lower bound of this functional
$E_1$ (without the knowledge of K\"ahler Einstein metric or
constant scalar curvature metric.).   Unfortunately,  even in the
simplest case of $CP^n (n> 1)$, we even do not know how to prove
directly why $E_1$ energy is bounded from below
 without using Mabuchi-Bando theorem.  On $S^2$,  the previous proof (for example, given by
 Osgood-Phillips-Sarnak  \cite{OPS})
 relies on the existence of constant scalar curvature metric (as we showed before).  The present paper
gives an aprori lower bound of $E_1$ (or Liouville energy)  on the
topological $S^2.\;$ This certainly sheds light on how to
analytically derive the low bound of the $E_1$ in higher dimension
and it   gives us hope that the large program can be successful.

\medskip
\subsection{ The main results}
The main purpose of this note is to give such an analytic proof of

\begin{thm} \label{main}
If  $(M, g)$ is a topological two sphere with volume $4\pi$ and
bounded curvature, then there is a  constant $C_{low}$ (which
depends on metric $g$) such that for any conformal metric
$g_1=e^{2u} g$ with volume $4\pi$,
$$
L(g_1) \ge C_{low}.
$$
\end{thm}

\begin{rem}The previous proof  based on the uniformization theorem and the Onofri inequality
 does yield   the best constant for
$C_{low}$.  Our proof can not yield the best lower bound.
\end{rem}


Comparing with some existing results for higher dimensional Yamabe
type problem (see, for example, \cite{HV} and \cite{LZ1}), our
proof seems quite standard. We start with a local sharp inequality
(caution: this local inequality does not rely on Moser's
inequality), which allows us to obtain a {\it rough} sharp
inequality (or, sometimes, we refer it as  an $\E$-level
inequality) in Proposition \ref{prop2-1}. We then study the
behavior  of the minimizing sequence as the parameter $\E$ goes to
$0$. It will be shown that the perturbed Liouville energy will
stay uniformly bounded as $\E$ goes to $0$, and this yileds the
result in Theorem \ref{main}. It worths of mentioning that in
higher dimensional case, similar argument was carried out by Hebey
and Vaugon \cite{HV} in their proof of a conjecture due to Aubin
\cite{A}, and by Y.Y.Li and M.Zhu in the proof of sharp trace
inequality \cite{LZ1}.

\bigskip

\section {Preliminary results}

We start with a local sharp inequality which is proved in
\cite{LZ}.  It is interesting to point out that such a local
inequality also yields a much simpler proof of the Onofri
inequality (see details in \cite{LZ}).

Let $\Omega \subset R^2$ be a bounded smooth domain, $\Omega^*$ be
the ball in $R^2$ which has the same area as $\Omega$, and denote
$$
D_{a,b}(\Omega)=\{f(x)-b\in H^1_0(\Omega) \ : \ \int_{\Omega}
e^{2f}dx=a\}.
$$
We have the following sharp  inequality.
\begin{thm}\cite{LZ}
$$\inf_{w \in D_{a,b}(\Omega)}  \int_{\Omega}|\nabla w|^2 d x\ge 4
\pi \cdot (\ln \frac{a e^{-2b}}{\pi r^2} +\frac {\pi r^2}{a
e^{-2b}}-1),
$$
where $r$ is the radius of $\Omega^*$. \label{thm2-1}
\end{thm}
As simple consequences, we have
\begin{cor} Let $\M \subset R^2$ be a smooth bounded domain.  For any $u\in H_0^1(\Omega)$,
$$
\int_{\Omega} e^{u}dx \le \pi r^2 e\exp\{\frac 1{16 \pi}
  \int_{\Omega}|\nabla u|^2 d x\},
$$
where  $r$ is the radius of $\Omega^*$.\label{cor2-1}
\end{cor}
And
\begin{cor} Let  $(M, g)$ be a smooth  Riemann surface. For any given
$\epsilon>0$, there are $r_0>0$  and $C(\E)>0$ such that if
$\Omega \subset B(x_0, r_0)$ for some point $x_0\in M$, then
$$
\int_{\Omega} e^{w}dV_g \le C(\E)\exp\{(\frac 1{16 \pi}+\epsilon)
  \int_{\Omega}|\nabla_g w|^2 d V_g\} \ \ \ \ \ \forall  w\in H_0^1(\Omega).
$$
\label{cor2-2}
\end{cor}

\begin{rem}
Inequalities in Corollary \ref{cor2-1}-\ref{cor2-2} were proved
early by  Cherrier \cite{Ch} based on Moser's inequality.
\label{reme-1}
\end{rem}

From these two corollaries we can derive
\begin{lem} Let  $(M, g)$ be a given   Riemann surface. For any given
$\epsilon>0$, there are  constants $C_1=C_1(\epsilon), \
C_2=C_2(\epsilon)$ such that
$$
\int_{M} e^{w}dV_g \le C_1 \exp\{(\frac 1{16 \pi}+\epsilon)
  \int_{M}|\nabla_g w|^2 d V_g+C_2\int_M w^2\}, \ \ \ \ \ \forall  w\in H^1(M).
$$
\label{lem2-1}
\end{lem}

\noindent{\bf Proof.}\ Let $r_0$ be the constant given in
Corollary \ref{cor2-2}. Let $\bigcup_{i=1}^N B(r_i/2, x_i)$ be a
covering of $M$, where $r_i <r_0$. Let $\{\psi_i(x)\}_{i=1}^N$ be
 smooth cutoff functions subordinate to this covering,
satisfying:
$$
\psi_i(x)= \left\{
\begin{array}{rll}
&1, \ \ \ x \in B(x_i, r_i/2)\\
&0, \ \ \ x \notin B(x_i, r_i),
\end{array}
\right.
$$
and $0\le \psi_i^2\le 1.$ We have, using Corollary \ref{cor2-2},
that
$$
\begin{array}{rll}
\int_M e^u d V_g&\le \sum_{i=1}^N \int_{B(x_i, r_i/2)} e^u d V_g\\
& \le \sum_{i=1}^N \int_{B(x_i, r_i)} e^{u\cdot \psi_i}d V_g\\
&\le \sum_{i=1}^N  C(\epsilon)\exp\{(\frac 1{16 \pi}+\epsilon)
  \int_{ B(x_i, r_i)}|\nabla_g (u\cdot \psi_i)|^2 d V_g\}.
  \end{array}
$$
Note
$$
\begin{array}{rll}
 \int_{ B(x_i, r_i)}|\nabla_g (u\cdot \psi_i)|^2 d V_g \le &\int_{ B(x_i, r_i)}
 (|\nabla_g u|^2 \cdot |\psi_i|^2+C_2 u^2) d V_g\\
 \le & \int_{M}
 (|\nabla_g u|^2 +C_2 u^2 )d V_g.
 \end{array}
 $$
 It follows that
 $$
 \int_M e^u dV_g \le N\cdot  C(\epsilon) \exp\{(\frac 1{16 \pi}+\epsilon)
  \int_{M}  (|\nabla_g u|^2 +C_2 u^2) d V_g.
  $$

 \medskip

  We are now ready to establish the following
  $\epsilon-$level inequality according to Aubin \cite{A2}. Similar type inequalities were
  discussed for higher dimensional cases. See, for example,
  \cite{HV} or \cite{LZ1}.

\begin{prop} Let $(M, g)$ be a  topological two sphere with bounded Gaussian
curvature $K_g$. Then for any given $\epsilon>0$, there is a
constant $C_3=C_3(\epsilon)$ such that for any $u\in H^1(M)$ with
$\int_M K_g \cdot udV_g =0$, \bg \int_M e^u dV_g \le C_3(\epsilon)
\exp\{(\frac 1{16 \pi}+\epsilon)
  \int_{M}  |\nabla_g u|^2 d V_g\}.
\label{2-1-1} \ede \label{prop2-1}
\end{prop}

To prove this proposition, we need the following
Poincar\'e-Sobolev type inequality.

\begin{lem}   Let $(M, g)$ be a  topological two sphere with bounded Gaussian
curvature $K_g$. Then for any $p\ge 1$, there is a constant $c_p$
such that for any $u \in H^1(M)$ with  $\int_M K_g \cdot udV_g
=0$, \bg (\int_M |u|^p)^{2/p} \le c_p \int_M |\nabla_g u|^2 dV_g.
\label{2-1} \ede \label{lem2-2}
\end{lem}

\noindent{\bf Proof.}\ The proof is standard. We prove it by
contradiction. Suppose this is not true, then there is a sequence
of functions $\{u_n\}$ in $H^1(M)$ with  $\int_M K_g \cdot u_ndV_g
=0$ such that
$$
(\int_M |u_n|^p)^{2/p}\ge n \int_M |\nabla_g u_n|^2 dV_g.$$ Let $$
v_n=\frac {u_n}{(\int_M |u_n|^p)^{1/p}}.$$ We know that
$$
||v_n||_{L^p}=1, \ \ \ ||v_n||_{H^1}\le C, \ \ \mbox{and} \ \
||\nabla_g v_n||_{L^2} \to 0.$$ Thus, $v_n \rightharpoonup v_0$ in
$H^1(M)$ and $v_0$ is a constant. By the compact embedding, we
know that $v_n \to v_0$ in $L^q$ for any $q<\infty$, thus $\int_M
K_g v_0 dV_g =0$  and $||v_0||_{L^p}=1$. On the other hand, since
$\int_M K_g dV_g =4 \pi$, we conclude that $v_0=0$. Contradiction.

\begin{rem}
The topological assumption is obviously needed since the above
inequality does not hold for a flat torus.
\end{rem}

\begin{rem}The condition on $K_g$ could be relaxed to $
||K_g||_{L^q(M)}\le C$ for some $q>1$. However, it is not obvious
whether the result is still true or not by only assuming that $K_g
\in L^1(M)$.  \label{rem2-2}
\end{rem}

We now return to the proof of Proposition \ref{prop2-1}.

\medskip

\noindent{\bf Proof of Proposition \ref{prop2-1}.}\  Due to Lemma
1 in \cite{A} (on page 586), we only need to prove the proposition
for functions in $C^1(M)$ with only non-degenerate critical
points. Let $u(x) $ be such a function  with $\int_M K_g \cdot u
dV_g=0$. For small enough $\eta>0$, there exists $a_\eta$ such
that \bg Vol \{ x \in M \ : \ u(x)>a_\eta\} =\eta.\label{2-2} \ede
Using Lemma \ref{lem2-1} we have  for any $\epsilon_1>0$, there
are constants $C_1(\epsilon_1), \ C_2(\epsilon_1)>0$ such that \bg
\begin{array}{rll}
\int_M e^u dV_g  \le & e^{a_\eta} \cdot \int_M e^{(u-a_\eta)} dV_g
\le e^{a_\eta} \cdot
\int_M e^{(u-a_\eta)_+ }dV_g\\
\le & e^{a_\eta} \cdot C_1(\epsilon_1) \exp\{(\frac 1{16
\pi}+\epsilon_1)
  \int_{M}|\nabla_g (u-a_\eta)_+|^2 d V_g\\
  &+C_2(\epsilon_1)\int_M |(u-a_\eta)_+|^2 dV_g\},
\end{array}
\label{2-3} \ede where $f_+(x)=\max \{ f(x),  0\}.$ Since
$$ a_\eta \cdot \eta \le \int_{\{x\in M \ : u(x)>a_\eta\}} udV_g \le
\int_M|u| dV_g \le c_1^{1/2}(\int_M|\nabla_g u|^2)^{1/2},
$$
we have
$$
a_\eta \le \eta ||\nabla_g u||_{L^2}^2 +\frac C{\eta^3}.
$$
Also, using Lemma \ref{lem2-2}, we have
$$
\begin{array}{rll}
\int_M|(u-a_\eta)_+|^2 dV_g = & \int_{\{x\in M \ :
u(x)>a_\eta\}}|(u-a_\eta)_+|^2 dV_g \\
\le & (\int_{\{x\in M \ : u(x)>a_\eta\}}|(u-a_\eta)_+|^4
dV_g)^{1/2} \cdot \eta^{1/2} \\
\le &  (\int_{M}|u|^4
dV_g)^{1/2} \cdot \eta^{1/2}\\
\le& c_4||\nabla_g u||_{L^2}^2 \cdot \eta^{1/2}.
\end{array}
$$
If we choose $\epsilon_1 $ and $\eta$ sufficiently small so that
$$
\eta+\epsilon_1+C_2(\E_1) \cdot c_4 \eta^{1/2} <\epsilon, $$ we
obtain  inequality (\ref{2-1-1}).

\bigskip

\section{Proof of Theorem \ref{main}}

In this section, we always assume that $(M, g)$ is a smooth
topological sphere  with bounded Gaussian curvature $K_g$. For
convenience, we also assume that $\int_M dV_g =4 \pi$ and  use the
simplified  notations as follows:
$$ \nabla:=\nabla_g, \ \ \ \Delta:=\Delta_g, \ \ \ \int_M f :=\int_M f dV_g.$$ For any
$\epsilon>0$, we define a perturbed (modified also) Liouville
energy functional by \bg I_\epsilon (u)=\frac 12 \int_M |\nabla
u|^2 +\frac{8\pi -\epsilon}{ 4 \pi}\int_M K_g u-( 8\pi -\E)\ln
\int_M e^u, \label{3-1} \ede and denote
$$
E_\epsilon=\inf_{u \in H^1(M), \ \int_M K_g u =0} I_\E (u).
$$
It is easy to see that $I_\E(u+C)=I_\E(u)$ for any constant $C$.
Due to Proposition \ref{prop2-1} we know that $E_\E$ is a finite
number. More over, we can show

\begin{lem} $E_\E$ is achieved by a function $u_\E \in H^1(M)$.
\label{lem3-1}
\end{lem}

\noindent{\bf Proof.}\ Let $\{u_n\}$ be a minimizing sequence with
$\int_M K_g u_n=0$. Choose a small positive number $\E_1$ so that
$$ (\frac 12-\E_1)\cdot \frac 1{8 \pi-\E}>\frac 1{16 \pi}+\E_1.$$
 Using Proposition \ref{prop2-1} we have
$$
\begin{array}{rll}
E_\E \leftarrow I_\E(u_n)&=\frac 12 \int_M |\nabla u_n|^2 -( 8\pi -\E)\ln \int_M e^{u_n}\\
&=(\frac 12-\E_1) \int_M |\nabla u_n|^2 -( 8\pi
-\E)\ln \int_M e^u+\E_1 \int_M |\nabla u_n|^2 \\
&\ge C(\E_1)+\E_1 \int_M |\nabla u_n|^2,
\end{array}
$$
thus $||\nabla u_n||_{L^2}<C.$ It follows from Lemma \ref{lem2-2}
that $||u_n||_{H^1}<C.$ Further, it implies that $u_n
\rightharpoonup u_\E$ in $H^1(M)$ and $e^{ku_n} \rightarrow
e^{ku_\E}$ in $L^1(M)$ for any positive number $k$ by Trudinger's
inequality (see, e.g. Theorem A in \cite{LZ}). It follows that
$u_\E$ is the minimizer.

\bigskip

Let $$ v_\E= u_\E-\ln\int_M e^{u_\E},$$
 then $v_\E$ is also a
minimizer, which satisfies the following Euler-Lagrange equation
\bg -\Delta v_\E+\frac{8 \pi-\E}{4 \pi}\cdot K_g =(8
\pi-\E)e^{v_\E} \ \ \ in \ \ M, \label{3-3} \ede and $\int_M
e^{2v_\E} \le C$ (since $||v_\E||_{H^1} \le C$). The standard
elliptic estimates yield that $v_\E \in C^\alpha(M)$ for some
$\alpha \in (0, 1).$

To prove Theorem \ref{main}, we only need to show that $E_\E$ is
uniformly bounded from below as $\E\to 0$. We follow  closely the
similar arguments given early in \cite{DJLW}.

We argue by contradiction. Suppose that $E_\E$ is not bounded from
below, that is, up to a subsequence (due to the nature of the
proof, for convenience we will not distinguish subsequence
$\{\E_i\}$ and the original sequence $\{\E\}$  in this section),
 \bg \lim_{\E \to 0} E_\E = -\infty.
 \label{3-4}\ede
 Let $\B_\E=v_\E(x_\E):=\max_{x \in M} v_\E(x).$ We first claim
 \begin{lem}
$$\lim_{\E \to 0} \B_\E=+\infty.$$
 \label{lem3-2}
 \end{lem}

\noindent {\bf Proof.}\ Let $v^a_\E=\frac 1{ 4\pi}\int_M v_\E$. If
$\B_\E<C$, then
 $e^{v_\E}<C$, thus $||\nabla (v_\E- v^a_\E)||_{L^2}<C$ from equation (\ref{3-3}).
 From Poincar\'e inequality, we know that $||v_\E- v^a_\E||_{L^2}<C$.  Since $||\Delta(v_\E- v^a_\E)||_{L^\infty}
 \le C$,
 we know from the standard elliptic estimates that $v_\E- v^a_\E
 \in C^1(M).$ If  $v^a_\E \to -\infty$ is unbounded, we know from equation (\ref{3-3}) that
  $ v_\E- v^a_\E \to w$, where $w$ satisfies $\Delta w=2K_g,$ this is impossible since $\int_M K_g \ne 0$;
  If $v^a_\E$ is bounded, then $v_\E=(v_\E-v^a_\E)+v^a_\E$ is bounded, thus
 $E_\E$ is bounded from below, this  is a contradiction to (\ref{3-4}).

 \medskip

Assume that, up to a subsequence, $x_\E \to \bar x\in M$.
  In a neighborhood of  $\bar x$, we choose
a normal coordinate system, and define
$$
\varphi_\E(x)=v_\E(\tau_\E^{-1}x+x_\E)-2\ln \tau_\E,$$ where
$\tau_\E:=e^{\B_\E/2},$ and without confusion, we use
$\tau_\E^{-1}x+x_\E$ to represent $exp_{ x_\E}(\tau^{-1}_\E x)$.
For any given $R>0$, if $\E$ is sufficiently small, $\varphi_\E$
satisfies \bg -\Delta \varphi_\E+\frac{8 \pi-\E}{4 \pi
\tau_\E^2}\cdot K_g =(8 \pi-\E)e^{\varphi_\E} \ \ \ in \ \
B_{2R}(0) \subset R^2. \label{3-5} \ede The behavior of $v_\E$ in
a tiny neighborhood of $\bar x$ can be described by the behavior
of $\varphi_\E$ in a large ball $B_{2R}(0) \subset R^2$.

\begin{lem} For a fixed $R>0$, there is a constant $C(R)$ such that
$$
|\varphi_\E(x)| <C(R) \ \ \ \ \forall \ \ x \in B_{R}(0).
$$
\label{lem3-3}
\end{lem}

\noindent{\bf Proof.}\ Let $\varphi_\E^{(1)}$ be the unique
solution to
$$
\left\{
\begin{array}{rll}
-\Delta \varphi_\E^{(1)} +\frac{8 \pi-\E}{4 \pi \tau_\E^2}\cdot
K_g &=(8 \pi-\E)e^{\varphi_\E} \ \ \ & in \ \ B_{2R}(0)\\
\varphi_\E^{(1)} &=0   \ \ \ & on \ \ \partial B_{2R}(0).
\end{array}
\right.
$$
Since $e^{\varphi_\E }\le 1$, we know that
$|\varphi_\E^{(1)}|<C_1(R)$ for $x \in B_R(0)$.  Let
$\varphi_\E^{(2)}= \varphi_\E- \varphi_\E^{(1)}$. Then
$\varphi_\E^{(2)}$ is a harmonic function bounded above by
$C_1(R)$. Applying Harnack inequality to
$2C_1(R)+\varphi_\E^{(2)}$ we have
$$
\frac 1{C} \le \frac
{2C_1(R)+\varphi_\E^{(2)}(x)}{2C_1(R)+\varphi_\E^{(2)}(0) }<C \ \
\ \ \ \forall x \in B_{R}(0).$$ This and the bound of
$\varphi_\E^{(1)}$  yield Lemma \ref{lem3-3}.

\bigskip

Since $\varphi_\E$ is uniformly bounded in any fixed ball, we
know, based on the standard elliptic estimates, that \bg
\varphi_\E \to \varphi_0 \ \ \ \ \mbox{in} \
C^{1,\alpha}(B_{R}(0)),\label{e-1} \ede for some $\alpha \in (0,
1)$,  where $\varphi_0$ satisfies \bg -\Delta \varphi_0=8 \pi
e^{\varphi_0} \ \ \  \mbox{in} \ \ R^2. \label{3-6} \ede On the
other hand, we note that \bg \varphi_0(x)\le \varphi_0(0)=1 \ \
\mbox{and} \ \ \int_{R^2} e^{\varphi_0} \le \lim_{R\to +\infty}
(\overline {\lim_{\E \to 0}} \int_{B_R} e^{\varphi_\E}) \le 1.
\label{3-7} \ede From the classification of constant curvature
metric on $S^2$ or the PDE's result of Chen and Li \cite{CL} we
know that \bg \varphi_0(x)=2\ln \frac{1}{1+\pi|x|^2}. \label{3-9}
\ede

Away from the singular point $\bar x$,  we have the following
global estimate.

\begin{lem}
For any compact domain $K \subset \subset M\setminus\{\bar x\}$,
there is a constant $C(K)$ such that $$||v_\E -\frac 1{4\pi}
\int_M v_\E||_{L^\infty(K)} <C(K). $$ \label{lem3-4}
\end{lem}

The proof of Lemma \ref{lem3-4} is based on a lemma due to Brezis
and Merle \cite{BM}. Such a lemma plays the same role as the sharp
Sobolev inequalities play in higher dimensional cases. To keep our
argument intact, we relegate the proof of this lemma to the end of
this section.

\medskip

Note that  for any $1\le q<2$, \bg ||\nabla
v_\E||_{L^q(M)}<C_q.\label{e-2} \ede This is because for any
$\varphi \in W^{1,q/(q-1)}(M)$ with  $\int_M \varphi =0$ and
$||\varphi ||_{W^{1, q/(q-1)}}=1$ (thus $\varphi \in
L^\infty(M)$),
$$ |\int_M \nabla v_\E \nabla \varphi|=|\int_M \Delta v_\E
\varphi| \le C.$$ We derive from the above and (\ref{3-4}) that
\bg
\begin{array}{rll}
4 \pi v_\E^a = & \int_M K_g v_\E- \int_M K_g(v_\E-v_\E^a) \\
\le & \int_M K_g v_\E -C \to -\infty. \end{array}\label{e-3}\ede

From Lemma \ref{lem3-4} and (\ref{e-3}), we conclude, via the
standard elliptic estimates, that \bg v_\E - v_\E^a \to G(x, \bar
x) \ \ \ in \ \ C^{1,\alpha}(K) \label{3-10} \ede for any compact
domain $K$ in $M\setminus\{\bar x\},$ where $G(x, \bar x)$ is the
Green's function in $M$ satisfying \bg \left\{
\begin{array}{rll} -\Delta
G+2K_g&=8\pi \delta_{\bar
x} \ \ \ \ \ \mbox{in} \ \ \ M\\
\int_M G&=0.
\end{array}
\right.
 \label{3-12} \ede
Since $K_g$ is bounded, in a normal coordinate system centered at
$\bar x$, \bg G(x, \bar x)=-4\ln |x|+A(\bar x)+o(1), \label{3-13}
\ede where $A(\bar x)$ is a constant depending on the location of
$\bar x$, $o(1) \to 0$ as $|x| \to 0$.

To complete the proof of the main theorem, we also need a lower
bound for $v_\E$ away from the singular point $\bar x$.

\begin{lem} For any fixed $R>0$, let $r_\E=\tau_\E^{-1}R$. Then \bg v_\E(x) \ge G(x, \bar
x)-\B_\E-2 \ln \frac{1+\pi R^2}{R^2} -A(\bar x)+o_\E(1)
 \ \ \ \ \ \ \  \forall x
\in M\setminus B(\bar x, r_\E), \label{3-16} \ede where $o_\E(1)
\to 0$ as $\E \to 0$. \label{lem3-5}
\end{lem}

\noindent{\bf Proof.}\ On $\partial B(\bar x, r_\E)$, we define
$C^*_\E:=(v_\E-G) |_{|x|=r_\E}$. From (\ref{e-1}) and (\ref{3-9})
we know that
$$C^*_\E=-\B_\E-2 \ln \frac{1+\pi R^2}{R^2} -A(\bar
x)+o_\E(1).
$$
Let $$ \tilde K_g(x)=\left\{
\begin{array}{rll}
&K_g(x), & \mbox {if} \ \ K_g(x)\le 0\\
&\phi(x) K_g(x),  & \mbox{if} \ \ K_g(x)>0,
\end{array}
\right.
$$
where $0\le \phi \le 1$ is a measurable function such that $\int_M
\tilde K_g =0.$ Let $h(x)$ be one solution to $\Delta h=\tilde
K_g$ in $M$ (which is bounded).
 We  consider $v_\E(x) -G(x, \bar x)-C^*_\E+ \frac{\E}{4 \pi} h(x)$ in $M\setminus B(\bar
x, r_\E)$. Since $-\Delta (v_\E(x) -G(x, \bar x)-C^*_\E+
\frac{\E}{4 \pi} h(x)) \ge \frac {\E (K_g-\tilde K_g)}{4 \pi} \ge
0$,  we have  (\ref{3-16}) via the maximum principle.

\medskip

We are now ready to complete the proof of  the main theorem.

\medskip

We need to estimate $E_\E=I_\E(v_\epsilon)$. For any fixed  small
$\delta>0$, we assume that $\E$ is sufficiently small so that
$\delta>r_\epsilon$. Then
$$
\begin{array}{rll}
\int_M|\nabla v_\epsilon|^2&=\int_{M\setminus B(\bar x,
\delta)}|\nabla v_\epsilon|^2 +\int_{B(\bar x, \delta)\setminus
B(\bar x, r_\epsilon)} |\nabla v_\epsilon|^2+\int_{B(\bar x,
r_\epsilon)} |\nabla v_\epsilon|^2 \\
&:= {\mathbf I}_1+{\mathbf I}_2+{\mathbf I}_3.
\end{array}
$$

First of all, it is easy to see from (\ref{e-1}) and (\ref{3-9})
that \bg
\begin{array}{rll}
{\mathbf I}_3&=\int_{B(\bar x, r_\epsilon)}|\nabla
\varphi_0|^2dx+o_\epsilon(1)= \int_{B(R)}\left|\nabla
\left(2\ln\frac{1}{1+\pi|x|^{2}}\right)
\right|^2dx+o_\epsilon(1)\\
=&32\pi^3\int_0^R\frac{r^{2}r}{(1+\pi r^{2})^2}dr+o_\epsilon(1)
=16\pi\left( \ln(1+\pi R^{2})- 1 \right) +o_\epsilon(1)+o_R(1),
\end{array}
\label{3-17} \ede
 where $o_R(1) \to 0$ as $R\to \infty$.

 Since $v_\epsilon(x)-v^a_\E\to
G(x,\bar x)$ in $C^{2}_{loc}(M\setminus \{\bar x\})$, we obtain
\bg \begin{array}{rll}
 {\mathbf I}_1= \int_{M\setminus B(\bar x,
\delta)}|\nabla v_\epsilon|^2&=\int_{M\setminus B(\bar x,
\delta)}|\nabla G|^2+o_\epsilon(1) \\
&= - \int_{M\setminus B(\bar x, \delta)} 2 K_g G -\int_{\partial
B(\bar x, \delta)} G\frac{\partial G}{\partial{\mathbf
n}}+o_\epsilon(1).
\end{array}
\label{3-18} \ede

To estimate ${\mathbf I}_2$, we first use equation (\ref{3-3}) to
get
$$
\begin{array}{rll}
{\mathbf I}_2=&-\int_{ B(\bar x, \delta)\setminus B(\bar x,
r_\epsilon)}\Delta v_\epsilon\cdot v_\epsilon+
\int_{\partial(B(\bar x, \delta)\setminus B(\bar x, r_\epsilon))}
v_\epsilon
\frac{\partial v_\epsilon}{\partial{\mathbf n}} \\
&=\int_{ B(\bar x, \delta)\setminus B(\bar x, r_\epsilon)}(-\frac
{8 \pi-\E}{4\pi} K_g \cdot v_\epsilon+(8 \pi-\E)v_\epsilon
e^{v_\epsilon})+ \int_{\partial(B(\bar x, \delta)\setminus B(\bar
x, r_\epsilon))} v_\epsilon \frac{\partial
v_\epsilon}{\partial{\mathbf n}}.
\end{array}
$$
Applying Lemma \ref{lem3-5}, we have
$$
\begin{array}{rll}
 \int_{ B(\bar x, \delta)\setminus B(\bar x, r_\epsilon)}(8 \pi-\E)v_\epsilon
e^{v_\epsilon}& \ge -\B_\E \cdot (8 \pi-\E)\int_{ B(\bar x,
\delta)\setminus B(\bar x, r_\epsilon)}e^{v_\epsilon}\\
&+\int_{ B(\bar x, \delta)\setminus B(\bar x, r_\epsilon)}(8
\pi-\E)G e^{v_\epsilon}+o_\E(1)+o_R(1),
\end{array}
$$
where  $o_\E(1)\to 0$ as $\E \to 0$ and  $o_R(1)\to 0$ as $R \to
+\infty$. From equation (\ref{3-3}), (\ref{3-10}) and (\ref{3-12})
we obtain
$$
\begin{array}{rll}
\int_{ B(\bar x, \delta)\setminus B(\bar x, r_\epsilon)}(8
\pi-\E)G e^{v_\epsilon} =&-\int_{ B(\bar x, \delta)\setminus
B(\bar x, r_\epsilon)}2K_g {v_\epsilon} + \int_{ B(\bar x,
\delta)\setminus B(\bar x, r_\epsilon)}\frac {8 \pi-\E}{4 \pi}
K_g\cdot G\\
&+\int_{\partial(B(\bar x, \delta)\setminus B(\bar x,
r_\epsilon))}(v_\epsilon \frac{\partial G}{\partial{\mathbf n}}- G
\frac{\partial v_\epsilon}{\partial{\mathbf n}})\\
=&-\int_{ B(\bar x, \delta)\setminus B(\bar x, r_\epsilon)}2K_g
{v_\epsilon}+\int_{\partial B(\bar x, \delta)}v^a_\epsilon
\frac{\partial G}{\partial{\mathbf n}}\\
& -\int_{\partial B(\bar x, r_\epsilon)}(v_\epsilon \frac{\partial
G}{\partial{\mathbf n}}- G \frac{\partial
v_\epsilon}{\partial{\mathbf n}})\\
&+o_\E(1)+o_\delta(1),
\end{array}
$$
where $o_\delta(1)\to 0$ as $\delta \to 0$,  and $$
\begin{array}{rll}
 -\B_\E \cdot (8 \pi-\E)\int_{ B(\bar x,
\delta)\setminus B(\bar x, r_\epsilon)}e^{v_\epsilon}=&-\B_\E
\cdot \int_{ B(\bar x, \delta)\setminus B(\bar x,
r_\epsilon)}\frac{8\pi -\E}{4
\pi} \cdot K_g\\
&+\B_\E \cdot \int_{\partial(B(\bar x, \delta)\setminus B(\bar x,
r_\epsilon))}\frac{\partial v_\E}{\partial{\mathbf n}}.
\end{array}
$$
Thus
$$
\begin{array}{rll}
{\mathbf I}_2\ge &-\frac {16\pi-\E}{4 \pi} \cdot \int_{ B(\bar x,
\delta)\setminus B(\bar x, r_\epsilon)}K_g {v_\epsilon}
+\int_{\partial(B(\bar x, \delta)\setminus B(\bar x, r_\epsilon))}
v_\E \frac{\partial v_\E}{\partial{\mathbf n}}\\
&+\int_{\partial B(\bar x, \delta)}v^a_\epsilon \frac{\partial
G}{\partial{\mathbf n}} -\int_{\partial B(\bar x,
r_\epsilon)}(v_\epsilon \frac{\partial G}{\partial{\mathbf n}}- G
\frac{\partial
v_\epsilon}{\partial{\mathbf n}})\\
&-\B_\E \cdot \int_{ B(\bar x, \delta)\setminus B(\bar x,
r_\epsilon)}\frac{8\pi -\E}{4 \pi} \cdot K_g +\B_\E \cdot
\int_{\partial(B(\bar x, \delta)\setminus B(\bar x,
r_\epsilon))}\frac{\partial v_\E}{\partial{\mathbf n}}\\
&+o_\E(1)+o_\delta(1)+o_R(1) \\
=&-\frac {16\pi-\E}{4 \pi} \cdot \int_{ B(\bar x, \delta)\setminus
B(\bar x, r_\epsilon)}K_g {v_\epsilon} -\B_\E \cdot \int_{ B(\bar
x, \delta)\setminus B(\bar x, r_\epsilon)}\frac{8\pi -\E}{4 \pi}
\cdot K_g \\
&+o_\E(1)+o_\delta(1)+o_R(1) \\
 &-\int_{\partial B(\bar
x, r_\epsilon)}  \frac{\partial v_\E}{\partial{\mathbf n}}\cdot
(v_\E-(G-\B_\E))-\int_{\partial B(\bar x, r_\epsilon)} v_\epsilon
\frac{\partial
G}{\partial{\mathbf n}}\\
&+  \int_{\partial B(\bar x, \delta)}v_\E \cdot\frac{\partial
v_\E}{\partial{\mathbf n}}+v^a_\E \cdot \int_{\partial B(\bar x,
\delta)}\frac{\partial G}{\partial{\mathbf n}} +\B_\E \cdot
\int_{\partial B(\bar x, \delta)}\frac{\partial
v_\E}{\partial{\mathbf n}}.
\end{array}
$$

We now estimate the boundary term in the right hand side of the
above inequality. From (\ref{e-1}), (\ref{3-9}) and Lemma
\ref{lem3-5}, we know that
$$
\begin{array}{rll}
-\int_{\partial B(\bar x, r_\epsilon)}  \frac{\partial
v_\E}{\partial{\mathbf n}}\cdot (v_\E-(G-\B_\E))\ge \frac{8 \pi^2
R^2}{1+\pi R^2 }\cdot (-A(\bar x)-2 \ln\frac{1+\pi R^2}{R^2})
+o_\E(1),
\end{array}
$$
and
$$
\begin{array}{rll}
 &-\int_{\partial B(\bar x,
r_\epsilon)} v_\epsilon \frac{\partial G}{\partial{\mathbf n}}\\
&=-\B_\E \int_{\partial B(\bar x, r_\epsilon)} \frac{\partial
G}{\partial{\mathbf n}}-16 \pi \ln(1+\pi R^2)+o_\E(1)+o_R(1)\\
&= \B_\E \int_{M\setminus B(\bar x, r_\epsilon)} \Delta
G -16 \pi \ln(1+\pi R^2)+o_\E(1)+o_R(1)\\
 &=2\B_\E \int_{M\setminus B(\bar x,
r_\epsilon)}K_g-16 \pi \ln(1+\pi R^2)+o_\E(1)+o_R(1).
\end{array}
$$
From equation (\ref{3-3}) we have
$$
\begin{array}{rll}
\B_\E \cdot \int_{\partial B(\bar x, \delta)}\frac{\partial
v_\E}{\partial{\mathbf n}}&=-\B_\E \cdot \int_{ M\setminus B(\bar
x, \delta)}
\Delta v_\E\\
&=-\frac{8\pi -\E}{4 \pi} \B_\E\int_{ M\setminus B(\bar x,
\delta)} K_g+\B_\E \int_{ M\setminus B(\bar x, \delta)}(8 \pi -\E)
e^{v_\E}\\
&\ge -\frac{8\pi -\E}{4 \pi} \B_\E\int_{ M\setminus B(\bar x,
\delta)} K_g.
\end{array}
$$
Similarly, we have
$$
\begin{array}{rll}
v^a_\E \cdot \int_{\partial B(\bar x, \delta)}\frac{\partial
G}{\partial{\mathbf n}} = -2 v^a_\E \int_{ M\setminus B(\bar x,
\delta)} K_g,
\end{array}
$$
and
$$
\begin{array}{rll}
v^a_\E \cdot \int_{\partial B(\bar x, \delta)}\frac{\partial
v_\E}{\partial{\mathbf n}}&= -\frac{8\pi -\E}{4 \pi} v^a_\E \int_{
M\setminus B(\bar x, \delta)} K_g+v_\E^a\int_{ M\setminus B(\bar
x, \delta)}(8 \pi -\E) e^{v_\E}\\
&=-\frac{8\pi -\E}{4 \pi} v^a_\E \int_{ M\setminus B(\bar x,
\delta)} K_g+ v_\E^a e^{v_\E^a}\int_{ M\setminus B(\bar x,
\delta)}(8 \pi -\E) e^{v_\E- v_\E^a}\\
&=-\frac{8\pi -\E}{4 \pi} v^a_\E \int_{ M\setminus B(\bar x,
\delta)} K_g+o_\E(1).
\end{array}
$$

For the interior terms, we note that
$$
\begin{array}{rll}
\int_{ B(\bar x, \delta)\setminus B(\bar x, r_\epsilon)}K_g
{v_\epsilon}&=v^a_\E \int_{ B(\bar x, \delta)\setminus B(\bar x,
r_\epsilon)}K_g +\int_{ B(\bar x, \delta)\setminus B(\bar x,
r_\epsilon)}K_g ({v_\epsilon}-v^a_\E)\\
&=v^a_\E \int_{ B(\bar x, \delta)\setminus B(\bar x,
r_\epsilon)}K_g +o_\delta(1),
\end{array}
$$
where we use the Poincar\'e inequality and (\ref{e-2}).

We conclude
\bg
\begin{array}{rll}
{\mathbf I}_2\ge &-(16\pi-\E)v^a_\E+ \frac{16\pi-\E}{4 \pi}\int_{
B(\bar x,
r_\E)}K_g \cdot v^a_\E \ \\&-8\pi A(\bar x)-16 \pi\ln \pi -16 \pi \ln(1+\pi R^2)\\
&+\E \B_\E \cdot (1-\frac 1{4 \pi}\int_{ B(\bar x,
r_\E)}K_g)\\
&+\int_{\partial  B(\bar x, \delta)} G\frac{\partial
G}{\partial{\mathbf n}}\\
 &+o_\E(1)+o_\delta(1)+o_R(1).
\end{array}
\label{3-19} \ede

We have, from (\ref{3-17})-(\ref{3-19}), that
$$
\begin{array}{rll}
{ I}_\E(v_\E)\ge &-4\pi A(\bar x)-8 \pi\ln \pi -8 \pi \\
&+\E \B_\E \cdot (\frac 1 2-\frac 1{8 \pi}\int_{ B(\bar x,
r_\E)}K_g)\\
&-\frac 12\E v^a_\E+ \frac{16\pi-\E}{8 \pi}\int_{ B(\bar x,
r_\E)}K_g \cdot v^a_\E \ \\
&  - \int_{M\setminus B(\bar x, \delta)}  K_g
G+\frac{8\pi-\E}{4 \pi} \int_M K_g \cdot(v_\E-v_\E^a)\\
&+o_\E(1)+o_\delta(1)+o_R(1).
\end{array}
$$

From  (\ref{e-3}), we know that $v^a_\E \to -\infty$. Also, since
$\B_\E \cdot  vol(B(\bar x, r_\E)) =o_\E(1)$, we know from Lemma
\ref{lem3-5}, that $
 v^a_\E \cdot vol(B(\bar x, r_\E)) =o_\E(1)$.  We finally have
 $$
  \lim_{\E \to 0} E_\E \ge -4\pi A(\bar x)-8 \pi\ln \pi -8 \pi+C.
 $$
for some constant $C$. Thus Theorem \ref{main} is proved.

\bigskip

We are left to prove   Lemma \ref{lem3-4}. We need  the following
lemma due to Brezis and Merle \cite{BM}. \begin{lem} Let $\Omega
\subset M$ be a bounded domain and $u$ be a solution to the
following equation:
$$ \left \{
\begin{array}{rll} -\Delta u & =f(x) &\ \mbox {in} \ \ \Omega\\
u&=0 & \ \mbox{on} \ \ \partial \Omega. \end{array} \right.
$$
If $f(x) \in L^1(\Omega)$, then for any $0<\delta<4 \pi$, there is
a constant $C(\delta)$ such that $$ \int_\Omega \exp\{\frac{(4
\pi-\delta)|u(x)|}{||f||_{L^1(\Omega)}}\}\le C(\delta).
$$\label{lem3-7}
\end{lem}

Lemma \ref{lem3-7} was originally proved for a bounded domain in
$R^2$, but it is easy to see that the same result holds on any
domain in a compact Riemann surface.  We now return to the

\noindent {\bf Proof of Lemma \ref{lem3-4}.}\ Let $K \subset
\subset M\setminus \{\bar x\}$.  We choose other two  compact sets
$K_1$ and $K_2$ in $M\setminus \{\bar x\}$ such that $K \subset
\subset K_1 \subset \subset K_2 \subset \subset M\setminus \{\bar
x\}$. It follows from (\ref{e-1}) that, \bg \lim_{\E \to 0}
\int_{K_2} e^{v_\E}=0. \label{3-20} \ede

Let $v_\E^1$ be the unique solution to \bg \left\{
\begin{array}{rll}
-\Delta v_\E^1 &=(8\pi -\E) e^{v_\E} \
\ &  \ \mbox{on} \ \ K_2\\
v_\E^1&=0 \ \ & \mbox{on} \ \ \partial K_2,
\end{array}
\right. \label{3-21} \ede
 and $v_\E^2= v_\E -v_\E^1-v_\E^a.$ From
(\ref{3-20}) and Lemma \ref{lem3-7} we know that $e^{|v_\E^1|} \in
L^p(K_2)$ for some $1<p<2$; thus $||v_\E^1||_{L^{p}(K_2)} \le C.$
 On the other hand, we note $v_\E^2$ satisfying $\Delta
v_\E^2=\frac{8 \pi -\E}{4 \pi} K_g$ in $K_2$, we have, from the
interior $L^p$ estimates, that for $K_1 \subset \subset K_2$,
$$
\begin{array}{rll}
||v_\E^2||_{L^\infty(K_1)} & \le C||v_\E^2||_{L^{p}(K_2)}\\
& \le C(||v_\E-v_\E^a||_{L^{p}K_2)}+||v_\E^1|_{L^{p}(K_2)})\\
&\le C(||\nabla v_\E||_{L^p(K_2)}+||v_\E^1|_{L^{p}(K_2)})\\
&\le C,
\end{array}
$$
where we also use (\ref{e-2}). In turn, we know that $\int_{K_1}
e^{pv_\E^1}\cdot e^{pv_\E^2}\le C$. Combining with (\ref{e-3}) we
have
$$
\begin{array}{rll}
\int_{K_1} e^{pv_\E} & \le \int_{K_1} e^{pv_\E^a} \cdot
e^{pv_\E^1}
\cdot e^{pv_\E^2}\\
& \le C,
\end{array}
$$
we then know from (\ref{3-21}) that $||v_\E^1||_{L^\infty(K)}\le
C$, thus $||v_\E-v_\E^a||_{L^\infty(K)}\le C.$

\end{document}